\documentclass[11pt,twoside,reqno]{amsart}
\usepackage{amsmath,amsfonts,amsthm,amssymb}
\usepackage{graphicx,psfrag,overpic}
\usepackage{bm}
\usepackage[all,cmtip,line]{xy}
\usepackage{array}
\usepackage{color}
\numberwithin{equation}{section}
\numberwithin{figure}{section}

\usepackage{graphics}
\usepackage{epsfig}
\newtheorem{claim}{\bf \t}[part]

\newtheorem{theorem}{Theorem}[section]

\newtheorem{proposition}[theorem]{Proposition}
\newtheorem{remark}{Remark}[section]

\def\t{\theta}

\newcommand \R{\mathbb{R}}

\begin{document}
\title[Incompressible Limit for M-D Steady Compressible Euler Equations]
{Incompressible Limit of Solutions of Multidimensional Steady Compressible Euler Equations}
\author{Gui-Qiang G. Chen}
\address{G.-Q. Chen,
	Academy of Mathematics and Systems Science,
	Academia Sinica, Beijing 100190, P. R. China;
    School of Mathematical Sciences, Fudan University,
	Shanghai 200433, P. R. China;
	Mathematical Institute, University of Oxford,
	Radcliffe Observatory Quarter, Woodstock Road, Oxford, OX2 6GG, UK}
\email{chengq@maths.ox.ac.uk}

\author{Feimin Huang}
\address{F. Huang,  Academy of Mathematics and Systems Science,
	Academia Sinica, Beijing 100190, P. R. China}
\email{fhuang@amt.ac.cn}

\author{Tian-Yi Wang}
\address{T.-Y. Wang, Department of Mathematics, School of Science, Wuhan University
of Technology, Wuhan, Hubei 430070, P. R. China; Gran Sasso Science Institute, viale Francesco Crispi, 7, 67100 L'Aquila, Italy; The Institute of Mathematical
Sciences, The Chinese University of Hong Kong, Shatin, N.T., Hong Kong; Academy
of Mathematics and Systems Science, Academia Sinica, Beijing 100190, P. R. China;
Mathematical Institute, University of Oxford, Radcliffe Observatory Quarter, Woodstock
Road, Oxford, OX2 6GG, UK}
\email{tianyiwang@whut.edu.cn; tian-yi.wang@gssi.infn.it; wangtianyi@amss.ac.cn}

\author{Wei Xiang}
\address{W. Xiang, Department of Mathematics, City University of Hong Kong, Kowloon,
	Hong Kong, P. R. China}
\email{weixiang@cityu.edu.hk}
\date{\today}

\begin{abstract}
	A compactness framework is formulated for the incompressible limit
	of approximate solutions with weak uniform bounds with respect to
    the adiabatic exponent
	for the steady Euler equations for compressible fluids
	in any dimension.
	One of our main observations is
	that the compactness can be achieved
	by using only natural weak estimates for the mass conservation and
	the vorticity.
	Another observation is that the incompressibility of the
	limit for the homentropic Euler flow is directly from the continuity equation,
	while the incompressibility of the limit for the full Euler flow is
	from a combination of all the Euler equations.
    As direct applications of the compactness framework,
    we establish two incompressible limit theorems for
	multidimensional steady Euler flows through infinitely long nozzles,
which lead to two new existence theorems for the corresponding problems
for multidimensional steady incompressible Euler equations.
\end{abstract}

\maketitle

%%%%%%%%%%%%%%%%%%%%%%%%%%%%%%%%%%%%%%%%%%%%%%%%%%%%%%%%%%%%%%%%%%%%%%%%%%%%%%%%%%%%%%%%%%%%%%%%%%
\section{Introduction}
We are concerned with the incompressible limit of solutions of multidimensional
steady compressible Euler equations.
The steady compressible full Euler equations take the form:
\begin{eqnarray}\label{CFE}
	\begin{cases}
		\mbox{div}\,(\rho u)=0,\\
		\mbox{div}\,(\rho u\otimes u)+ \nabla p =0,\\
		\mbox{div}\, (\rho u E+ u p)=0,
	\end{cases}
\end{eqnarray}
while the steady homentropic Euler equations have the form:
\begin{eqnarray}\label{CISE}
	\begin{cases}
		\mbox{div}\,(\rho u)=0,\\[1mm]
		\mbox{div}\,(\rho u\otimes u)+ \nabla p =0,
	\end{cases}
\end{eqnarray}
where $x:=(x_1, \cdots, x_n)\in \R^n$ with $n\ge 2$, $u:=(u_1,\cdots,u_n)\in \R^n$ is the flow velocity,
\begin{equation}\label{speed}
	|u|=\Big(\sum_{i=1}^n u_i^2\Big)^{1/2}
\end{equation}
is the flow speed,  $\rho$, $p$, and $E$ represent the density, pressure, and total energy respectively,
and $u\otimes u:=(u_iu_j)_{n\times n}$ is an $n\times n$ matrix.

For the full Euler case, the total energy is
\begin{equation}\label{total-energy}
	E=\frac{|u|^2}{2}+\frac{p}{(\gamma-1)\rho},
\end{equation}
with adiabatic exponent $\gamma>1$, the local sonic speed is
\begin{equation}\label{sonic-speed-2}
	c=\sqrt{\frac{\gamma p}{\rho}},
\end{equation}
and the Mach number is
\begin{equation}\label{March-number-2}
	M=\frac{|u|}{c}=\frac{1}{\sqrt{\gamma}}\, |u|\sqrt{\frac{\rho}{p}}.
\end{equation}

For the homentropic case, the pressure-density relation is
\begin{equation}\label{pressure-1}
	p=\rho^\gamma, \qquad \gamma>1.
\end{equation}
The local sonic speed is
\begin{equation}\label{sonic-speed-1}
	c=\sqrt{\gamma \rho^{\gamma-1}}=\sqrt{\gamma}\, p^{\frac{\gamma-1}{2\gamma}},
\end{equation}
and the Mach number is defined as
\begin{equation}\label{March-number-1}
	M=\frac{|u|}{c}=\frac{1}{\sqrt{\gamma}}\,|u|p^{\frac{1-\gamma}{2\gamma}}.
\end{equation}

The incompressible limit is one of the fundamental fluid dynamic limits in fluid mechanics.
Formally, the steady compressible full Euler equations \eqref{CFE}
converge to the steady inhomogeneous incompressible Euler equations:
	\begin{eqnarray}\label{ICIHE}
		\begin{cases}
			\mbox{div}\, u=0,\\
			\mbox{div}\,(\rho u)=0,\\
			\mbox{div}\,(\rho u\otimes u)+ \nabla p=0,
		\end{cases}
	\end{eqnarray}
	while the homentropic Euler equations \eqref{CISE}
	converge to the steady homogeneous incompressible Euler equations:
	\begin{eqnarray}\label{ICHE}
		\begin{cases}
			\mbox{div}\, u=0,\\[1mm]
			\mbox{div}\,(u\otimes u)+ \nabla p=0.
		\end{cases}
	\end{eqnarray}
However, the rigorous justification of this limit for weak solutions has been
a challenging mathematical problem, since it is a singular limit
for which singular phenomena usually occur in the limit process.
In particular, both the uniform estimates and the convergence
of the nonlinear terms in the incompressible models are usually difficult to obtain.
Moreover, tracing the boundary conditions of the solutions
in the limit process is a tricky problem.
	
Generally speaking, there are two processes for the incompressible limit:
The adiabatic exponent $\gamma$ tending to infinity, and the Mach number $M$
tending to zero \cite{Masmoudi1, Masmoudi2}.
The latter is also called the low Mach number limit.
A general framework for the low Mach number limit for local smooth solutions
for compressible flow was established in Klainerman-Majda \cite{KM1,KM2}.
In particular,
the incompressible limit
of local smooth solutions of the Euler equations
for compressible fluids was established with  well-prepared initial data
{\it i.e.},
the limiting velocity satisfies the incompressible condition initially,
in the whole space or torus.
Indeed, by analyzing the rescaled linear group
generated by the penalty operator ({\it cf.} \cite{Ukai}),
the low Mach number limit can also be verified for the case of general data,
for which the velocity in the incompressible fluid is the limit
of the Leray projection of the velocity in the compressible fluids.
This method also applies to global weak solutions of the isentropic
Navier-Stokes equations with general initial data and various boundary
conditions \cite{Desjardins,Desjardins1,Lions P.-L.01}.
In particular, in \cite{Lions P.-L.01},
the incompressible limit on the stationary Navier-Stokes equations
with the Dirichlet boundary condition was also shown,
in which the gradient estimate on the velocity played the major role.
For the one-dimensional Euler equations,
the low Mach number limit has been proved by using the $BV$ space in \cite{CCZ}.
For the limit $\gamma\rightarrow \infty$,
it was shown in \cite{Lions-Masmoudi} that
the compressible isentropic Navier-Stokes flow would converge to
the homogeneous incompressible Navier-Stokes flow.
Later, the similar limit from the Korteweg barotropic Navier-Stokes model
to the homogeneous incompressible Navier-Stokes model was also
considered in \cite{Labbe}.
	
For the steady flow, the uniqueness of weak solutions of the steady
incompressible Euler equations is still an open issue.
Thus, the incompressible limit of the steady Euler
equations becomes more fundamental mathematically;
it may serve as a selection principle of physical relevant
solutions for the steady incompressible Euler equations
since a weak solution should not be regarded as the compressible
perturbation of the steady incompressible Euler flow in general.
Furthermore, for the general domain, it is quite challenging to obtain
directly a uniform estimate for the Leray projection of the velocity in
the compressible fluids.

In this paper, we formulate a suitable compactness framework
for weak solutions with weak uniform bounds with respect to
the adiabatic exponent $\gamma$ by employing the weak convergence argument.
One of our main observations is that the compactness can be achieved
by using only natural weak estimates for the mass conservation and
the vorticity, which was introduced in \cite{ChenHuangWang,Huang-Wang-Wang}.
Another observation is that the incompressibility of the
limit for the homentropic Euler flow follows directly from the continuity equation,
while the incompressibility of the limit for the full Euler flow is
from a combination of all the Euler equations.
Finally, we find a suitable framework to satisfy the boundary condition
without the strong gradient estimates on the velocity.
As direct applications of the compactness framework,
we establish two incompressible limit theorems for
multidimensional steady Euler flows through infinitely long nozzles.
As a consequence, we can establish the new existence
theorems for the corresponding problems for multidimensional steady
incompressible Euler equations.
	
The rest of this paper is organized as follows.
In \S 2, we establish the compactness framework for the incompressible
limit of approximate solutions of the steady full Euler equations and
the homentropic Euler equations in $\R^n$ with $n\geq2$.
In \S 3, we give a direct application of the compactness framework
to the full Euler flow through infinitely long nozzles in $\R^2$.
In \S 4, the incompressible limit of homentropic Euler flows
in the three-dimensional infinitely long axisymmetric nozzle is established.

\section{Compactness Framework for Approximate Steady Euler Flows}
	
In this section, we establish the compensated compactness framework for
approximate solutions of the steady Euler equations in $\mathbb{R}^n$
with $n\ge 2$.
We first consider the homentropic case, that is,
the approximate solutions $(u^{(\gamma)}, p^{(\gamma)})$ satisfy
\begin{eqnarray}\label{3.1}
	\begin{cases}
		\mbox{div}\big(\rho^{(\gamma)} u^{(\gamma)}\big)=e_1(\gamma),\\[1mm]
		\mbox{div}\big(\rho^{(\gamma)} u^{(\gamma)}\otimes u^{(\gamma)}\big)+ \nabla p^{(\gamma)}=e_2(\gamma),
	\end{cases}
\end{eqnarray}
where  $e_1(\gamma)$ and $e_2(\gamma):=(e_{21}(\gamma), \cdots, e_{2n}(\gamma))^\top$
are sequences of distributional functions depending on the parameter $\gamma$.
	
\begin{remark}
The distributional functions $e_i(\gamma)$, $i=1,2$, here present possible error terms
from different types of approximation.
If $(u^{(\gamma)}, p^{(\gamma)})$ with $\rho^{(\gamma)}:=\big(p^{(\gamma)}\big)^{\frac{1}{\gamma}}$
are the exact  solutions of the steady Euler flows,
$e_i(\gamma), i=1,2$, are both equal to zero.
Moreover, the same remark is true for the full Euler case, where $e_i(\gamma), i=1,2,3$,
are the distributional functions as introduced in \eqref{3.2}.
\end{remark}
	
Let the sequences of functions $u^{(\gamma)}(x):=(u^{(\gamma)}_1, \cdots, u^{(\gamma)}_n)(x)$
and  $p^{(\gamma)}(x)$ be defined on an open bounded
subset $\Omega\subset \mathbb{R}^n$ such that the following
qualities:
\begin{eqnarray}
&& \rho^{(\gamma)}:=(p^{(\gamma)})^{\frac{1}{\gamma}},\quad
|u^{(\gamma)}|:= \sqrt{\sum_{i=1}^n (u_i^{(\gamma)})^2}, \quad
		c^{(\gamma)} := \sqrt{\gamma}\big(p^{(\gamma)}\big)^{\frac{\gamma-1}{2\gamma}},\quad
		M^{(\gamma)} := \frac{|u^{(\gamma)}|}{c^{(\gamma)}},\label{2.2-a}\\[3mm]
&&E^{(\gamma)} = \frac{|u^{(\gamma)}|^2}{2}+\frac{\big(p^{(\gamma)}\big)^{\frac{\gamma-1}{\gamma}}}{\gamma-1}.
	\end{eqnarray}
can be well defined.
Moreover, the following conditions hold:
	
\medskip
(A.1). $M^{(\gamma)}$ are uniformly bounded by $\bar{M}$;
	
	\vspace{1.5mm}
	(A.2). $|u^{(\gamma)}|^2$ and $p^{(\gamma)}\ge 0$ are uniformly bounded in $L^1_{loc}(\Omega)$;
	
	\vspace{1.5mm}
	
	(A.3). $e_1(\gamma)$ and $\mbox{curl}\  u^{(\gamma)}$ are in a compact set in $H_{loc}^{-1}(\Omega)$;
	
	\vspace{1.5mm}
		
(H). As $\gamma \rightarrow \infty$,
$$
\int_{\Omega} \ln\big(E^{(\gamma)}\big)\, {\rm d}x =o(\gamma).
$$

\medskip
\smallskip
\begin{remark}
In the limit $\gamma\to \infty$,
the energy sequence $E^{(\gamma)}$ may tend to zero.
Condition {\rm (H)} is designed to exclude the case that $E^{(\gamma)}$ exponentially decays
to zero as $\gamma\to \infty$.
In fact, in the two applications in \S $3$--\S $4$ below, both of the energy sequences $E^{(\gamma)}$
{go to zero with polynomial rate} so that condition {\rm (H)} is satisfied automatically.
It is noted that condition {\rm (H)} could be replaced equivalently by a pressure condition:
$$
\int_{\Omega} \ln\big(p^{(\gamma)}\big)\, {\rm d}x =o(\gamma) \qquad \mbox{as}\,\, \gamma \rightarrow \infty.
$$
Indeed, from {\rm (A.1)} and \eqref{2.2-a}, we have
\begin{equation}\label{e}
\frac{1}{\gamma-1}(p^{(\gamma)})^{1-\frac1\gamma}
\le E^{(\gamma)}=\frac{|u^{(\gamma)}|^2}{2}+\frac{\big(p^{(\gamma)}\big)^{\frac{\gamma-1}{\gamma}}}{\gamma-1}
\le \frac{(\gamma-1)\gamma\bar{M}^2+2}{2(\gamma-1)}(p^{(\gamma)})^{1-\frac1\gamma},
\end{equation}
which directly implies the equivalence of the two conditions.
\end{remark}
	
\begin{remark}
Conditions {\rm (A.1)}--{\rm (A.3)} are naturally satisfied in the applications in \S $3$--\S $4$ below.
\end{remark}

Then we have
\begin{theorem}[Compensated compactness framework for the homentropic Euler case] \label{thm2.1}
Let a sequence of functions $u^{(\gamma)}(x)=(u^{(\gamma)}_1, \cdots, u^{(\gamma)}_n)(x)$
and $p^{(\gamma)}(x)$ satisfy conditions {\rm (A.1)}--{\rm (A.3)} and {\rm (H)}.
Then there exists a subsequence (still denoted by) $( u^{(\gamma)}, p^{(\gamma)})(x)$
such that, when $\gamma\rightarrow \infty$,
\begin{equation}
\begin{array}{ll}
				u^{(\gamma)}(x)\rightarrow (\bar{u}_1, \cdots, \bar{u}_n)(x) &\quad\mbox{{\it a.e.} in}\,\,\, x \in \Omega,
				\\[2mm]
				p^{(\gamma)}(x)\rightharpoonup \bar{p}\quad&
				\quad\mbox{in bounded measure.}
\end{array}
\end{equation}
\end{theorem}
	
\smallskip		
\noindent\textbf{Proof}.
We divide the proof into four steps.
	
\medskip
1. From condition
(A.2),
we can see that $p^{(\gamma)}$ weakly converges
to $\bar{p}$ in measure as $\gamma\rightarrow\infty$.
	
\medskip
2.  Now we show that $\rho^{(\gamma)}=(p^{(\gamma)})^{\frac{1}{\gamma}}(x)\rightarrow 1$
{\it a.e.} in $x \in \Omega$ as $\gamma\rightarrow \infty$.
	
\smallskip
Since $\gamma\to \infty$, for given $q\ge 1$,
we may assume $\gamma>q$.
Then we find by Jensen's inequality that
\begin{equation}
\left(\frac{\int_{K}(p^{(\gamma)})^{\frac{q}{\gamma}}\, {\rm d}x}{|\Omega|}\right)^{\frac{\gamma}{q}}
\le \frac{\int_{\Omega}p^{(\gamma)}\, {\rm d}x}{|\Omega|},
\end{equation}
where $|\Omega|$ is the Lebesgue measure of $\Omega$.
Then, for $(p^{(\gamma)})^{\frac{1}{\gamma}}$, we have
\begin{equation}\label{J1}
\left(\int_{\Omega}(p^{(\gamma)})^{\frac{q}{\gamma}}\, {\rm d}x\right)^{\frac{1}{q}}
\le \left(\int_{\Omega} p^{(\gamma)}\, {\rm d}x\right)^{\frac{1}{\gamma}}|\Omega|^{\frac{1}{q}-\frac{1}{\gamma}}.
\end{equation}
	
On the other hand, since $\ln y$ is concave with respect to $y$, we have
\begin{equation}
\frac{1}{|\Omega|}\int_{\Omega}\ln((p^{(\gamma)})^{\frac1\gamma})\,{\rm d}x \leq
\ln (\frac{1}{|\Omega|}\int_{\Omega}(p^{(\gamma)})^{\frac1\gamma}\,{\rm d}x),
\end{equation}
which implies from the H\"{o}lder inequality that
\begin{equation}\label{J2}
\begin{split}
|\Omega|^{\frac{1}{q}}\exp\Big\{\frac{1}{\gamma|\Omega|}\int_{\Omega}\ln (p^{(\gamma)})\, {\rm d}x\Big\}
&\le |\Omega|^{\frac{1}{q}}\exp\Big\{\ln \big(\frac{1}{|\Omega|}\int_{\Omega}(p^{(\gamma)})^{\frac1\gamma}\, {\rm d}x\big)\Big\}\\
&=|\Omega|^{\frac{1}{q}-1} \int_{\Omega}(p^{(\gamma)})^{\frac1\gamma}\, {\rm d}x\\
&\le \left(\int_{\Omega}(p^{(\gamma)})^{\frac{q}{\gamma}}\, {\rm d}x\right)^{\frac{1}{q}}.
\end{split}
\end{equation}
Moreover, from \eqref{e}, we have
\begin{equation}
\ln\left(E^{(\gamma)}\right)
\leq\frac{\gamma-1}{\gamma}\ln p^{(\gamma)}-\ln\big(\frac{2(\gamma-1)}{(\gamma-1)\gamma\bar{M}^2+2}\big),
\end{equation}
which, together with $(\ref{J1})$ and $(\ref{J2})$, gives
\begin{equation}
|\Omega|^{\frac{1}{q}}\exp\Big\{\frac{\int_{\Omega}
\big(\ln (E^{(\gamma)})+\ln(\frac{2(\gamma-1)}{(\gamma-1)\gamma\bar{M}^2+2})\big)\, {\rm d}x}{(\gamma-1)|\Omega|}\Big\}
		\le \|(p^{(\gamma)})^{\frac{1}{\gamma}}\|_{L^q(\Omega)}
		\le \left(\int_{\Omega} p^{(\gamma)} \, {\rm d}x\right)^{\frac{1}{\gamma}}
		|\Omega|^{\frac{1}{q}-\frac{1}{\gamma}}.
\end{equation}
Note that both the left and right sides of the above inequality tend to $|\Omega|^{\frac1q}$ as $\gamma\to \infty$,
owing to condition (H).  Then we have
\begin{equation}\label{lim}
\lim_{\gamma\to \infty}\|\rho^{(\gamma)}\|_{L^q(\Omega)}=|\Omega|^{\frac1q},
\end{equation}
where $\rho^{(\gamma)}:=(p^{(\gamma)})^{\frac{1}{\gamma}}$.
In particular, taking $q=1$ and $q=2$ respectively, we have
\begin{equation}\label{lim1}
\lim_{\gamma\to \infty}\|\rho^{(\gamma)}\|_{L^2(\Omega)}=|\Omega|^{\frac12},\quad \lim_{\gamma\to \infty}\|\rho^{(\gamma)}\|_{L^1(\Omega)}=|\Omega|.
\end{equation}
This implies that $\rho^{(\gamma)}$ are uniformly bounded in $L^2(\Omega)$.
Then there exists a subsequence of $\rho^{(\gamma)}$ (still denoted by $\rho^{(\gamma)}$) such that $\rho^{(\gamma)}$ weakly converges to $\bar{\rho}$ in $L^2(\Omega)$.
By a simple computation, we obtain from \eqref{lim1} that
$$
\int_{\Omega}(\bar{\rho}-1)^2dx=\int_{\Omega}(\bar{\rho}^2-2\bar{\rho}+1)\, {\rm d}x
\le \lim_{\gamma\to \infty} \int_{\Omega}\big((\rho^{(\gamma)})^2-2\rho^{(\gamma)}+1\big)\, {\rm d}x=0.
$$
That is, $\rho^{(\gamma)}$ converges to $1$ {\it a.e.} in $x \in \Omega$, as $\gamma\rightarrow \infty$.
	
\smallskip
3. By the div-curl lemma of Murat \cite{Murat} and Tartar \cite{Tartar},
the Young measure representation theorem for a uniformly bounded sequence
	of functions in $L^p$ ({\it cf.} Tartar \cite{Tartar}; also see Ball \cite{Ball}),
	we use $(\ref{3.1})_1$ and (A.3) to obtain the following commutation identity:
	\begin{equation}\label{2.3}
		\sum\limits_{i=1}^n\langle \nu(\rho, u), \, u_i\rangle \langle \nu(\rho, u), \, \rho u_i\rangle
		=\langle \nu(\rho, u), \,\sum\limits_{i=1}^n \rho u_i^2\rangle,
	\end{equation}
where we have used that $\nu(\rho, u)$ is the associated Young measure (a probability measure)
for the sequence  $(\rho^{(\gamma)}, u^{(\gamma)})(x)$.
	
Then the main point in the
compensated compactness framework is to prove that
$\nu(\rho, u)$ is in fact a Dirac measure, which
in turn implies the compactness of the sequence
$(\rho^{(\gamma)}, u^{(\gamma)})(x)$.
On the other hand, from
$$
\lim_{\gamma\rightarrow\infty}\rho^{(\gamma)}(x)=1 \qquad {~~ a.e.}
$$
we see that
$$
\nu(\rho, u)=\delta_{1}(\rho)\otimes\nu(u),
$$
where $\delta_1(\rho)$ is the Delta mass concentrated at $\rho=1$.

\smallskip
4. We now show $\nu(u)$ is a Dirac measure.

Combining both sides of $(\ref{2.3})$ together, we have
\begin{equation}
\langle \nu(u^{(1)})\otimes\nu(u^{(2)}), \,\sum\limits_{i=1}^n u^{(1)}_i (u_i^{(1)}-u_i^{(2)})\rangle=0.
\end{equation}
Exchanging indices $(1)$ and $(2)$, we obtain the following symmetric commutation identity:
\begin{equation}\label{3.4}
		\langle \nu(u^{(1)})\otimes\nu(u^{(2)}), \,
		\sum_{i=1}^n(u_i^{(1)}-u_i^{(2)})^2\rangle=0,
\end{equation}
which immediately implies that,
$$
u^{(1)}=u^{(2)},
$$
{\it i.e.}, $\nu(u)$ concentrates on a single point.

If this would not be the case, we could suppose {that there are} two different points
$\acute{u}$ and $\grave{u}$ in the support of $\nu$.
Then $(\acute{u}, \acute{u})$, $(\acute{u}, \grave{u})$, $(\grave{u}, \acute{u})$,
and $(\grave{u}, \grave{u})$
would be in the support of $\nu\otimes\nu$,
which contradicts with $u^{(1)}=u^{(2)}$.

Therefore, the Young measure $\nu$ is a Dirac measure, which implies the strong convergence
of $u^{(\gamma)}$. This completes the proof.

\bigskip
	
For the full Euler case, we assume that the approximate solutions $(\rho^{(\gamma)}, u^{(\gamma)},p^{(\gamma)})$ satisfy
\begin{eqnarray}\label{3.2}
\begin{cases}
\mbox{div}\big(\rho^{(\gamma)} u^{(\gamma)}\big)=e_1(\gamma),\\[1mm]
\mbox{div}\big(\rho^{(\gamma)} u^{(\gamma)}\otimes u^{(\gamma)}\big)+ \nabla p^{(\gamma)}=e_2(\gamma),\\[1mm]
\mbox{div}\big(\rho^{(\gamma)} u^{(\gamma)} E^{(\gamma)}+ u^{(\gamma)} p^{(\gamma)}\big)=e_3(\gamma),
\end{cases}
\end{eqnarray}
where $e_1(\gamma)$, $e_2(\gamma)=(e_{21}(\gamma), \cdots, e_{2n}(\gamma))^\top$,
and $e_3(\gamma)$ are sequences of distributional functions depending on the parameter $\gamma$.
In this case, {the energy function is
$$
E^{(\gamma)}:=\frac{|u^{(\gamma)}|^2}{2}+\frac{ p^{(\gamma)}}{(\gamma-1)\rho^{(\gamma)}},
$$
and} the entropy function is
$$
S^{(\gamma)} := \frac{\rho^{(\gamma)}}{(p^{(\gamma)})^{\frac{1}{\gamma}}}\ge 0,
$$
so that condition {\rm (H)} for the homentropic case is replaced by
	
\medskip
(F.1). As $\gamma\to \infty$,
$$
\int_{\Omega} \ln(p^{(\gamma)})\, {\rm d}x =o(\gamma) \qquad \mbox{as}\,\, \gamma \rightarrow \infty;
$$
	
\vspace{1mm}
(F.2).  $S^{(\gamma)}$ converges to a bounded function $\overline{S}$ {\it a.e.} in $\Omega$ as $\gamma\rightarrow\infty$.
	
\medskip
\vspace{2mm}

\begin{remark}
Conditions {\rm (A.1)}--{\rm (A.3)} and {\rm (F.1)}--{\rm (F.2)} in the framework
are naturally satisfied in the applications for the full Euler case
in \S $3$ below.
\end{remark}

Similar to Theorems \ref{thm2.1}, we have

\begin{theorem}[Compensated  compactness framework for the full Euler case] \label{thm3.1}
Let a sequence of functions $\rho^{(\gamma)}(x)$,
$u^{(\gamma)}(x)=(u^{(\gamma)}_1, \cdots, u^{(\gamma)}_n)(x)$,
and $p^{(\gamma)}(x)$ satisfy conditions {\rm (A.1)}--{\rm (A.3)} and {\rm (F.1)}--{\rm (F.2)}.
Then there exists a subsequence (still denoted by) $(\rho^{(\gamma)}, u^{(\gamma)}, p^{(\gamma)})(x)$
such that, as $\gamma\rightarrow\infty$,
\begin{equation}
\begin{array}{ll}
p^{(\gamma)}(x)\rightharpoonup \bar{p}\quad	&\quad \mbox{in bounded measure,}\\[1mm]
\rho^{(\gamma)}(x)\rightarrow \bar{\rho}(x)\quad
				&\quad \mbox{{\it a.e.} in}\,\,\, x \in \Omega, \\[1mm]
u^{(\gamma)}(x)\rightarrow (\bar{u}_1, \cdots, \bar{u}_n)(x)
\qquad & \quad \mbox{{\it a.e.} in}\,\,\, x \in \{x\,:\, \bar{\rho}(x)>0, x\in \Omega\}.\nonumber
			\end{array}
		\end{equation}
	\end{theorem}
	
\noindent\textbf{Proof}. We follow the same arguments as in the homentropic case.

First, the weak convergence of $p^{(\gamma)}$ is obvious.
On the other hand, we observe
that \eqref{J1} and \eqref{J2} still hold for the full Euler case.
Then, for any $\gamma>q\ge 1$,
\begin{equation}\label{J2-1}
|\Omega|^{\frac{1}{q}}\exp\Big\{\frac{1}{\gamma|\Omega|}\int_{\Omega}\ln (p^{(\gamma)})\, {\rm d}x\Big\}
\le \left(\int_{\Omega}(p^{(\gamma)})^{\frac{q}{\gamma}}\, {\rm d}x\right)^{\frac{1}{q}}
\le \left(\int_{\Omega} p^{(\gamma)} \, {\rm d}x\right)^{\frac{1}{\gamma}}
|\Omega|^{\frac{1}{q}-\frac{1}{\gamma}}.
\end{equation}

Thanks to condition (F.1), we obtain
\begin{equation}\label{lim}
\lim_{\gamma\to \infty}\|(p^{(\gamma)})^{\frac1\gamma}\|_{L^q(\Omega)}=|\Omega|^{\frac1q}.
\end{equation}
Taking $q=1$ and $q=2$ respectively and following the same line of argument as in the homentropic case,
we conclude that $(p^{(\gamma)})^{\frac1\gamma}$ converges to $1$ {\it a.e.} in $x \in \Omega$
as $\gamma\rightarrow \infty$.
Then, from condition (F.2),
$\rho^{(\gamma)}=S^{(\gamma)}(p^{(\gamma)})^{\frac{1}{\gamma}}$ converges
to $\bar{\rho}:=\bar{S}\ge 0$ {\it a.e.} in $x \in \Omega$.

The remaining proof is the same as that for the homentropic case,
except the strong convergence of $u^{(\gamma)}$ only stands on $ \{x\,:\, \bar{\rho}(x)>0, x\in \Omega\}$
since the vacuum can not excluded.
This completes the proof.

\medskip
\begin{remark}\label{rem4}
			Consider any function
			$Q(\rho, u, p):=(Q_1,\cdots,Q_n)(\rho, u, p)$ satisfying
			\begin{equation}\label{3.12}
			\mbox{\rm div}\,(Q(\rho^{(\gamma)}, u^{(\gamma)},  p^{(\gamma)}))=e_Q(\gamma),
			\end{equation}
			where $e_Q(\gamma)\rightarrow 0$ in the distributional sense as
			$\gamma\rightarrow \infty$.
			The similar statement is also valid for Theorem $\ref{thm3.1}$, via replacing  $(\ref{3.1})$
			by \eqref{3.2}.
\end{remark}

\medskip
Then, as direct corollaries, we conclude the following propositions.
	
\begin{proposition}[Convergence of approximate solutions of the homentropic Euler equations]\label{thm2.2}
Let $u^{(\gamma)}(x)=(u^{(\gamma)}_1, \cdots, u^{(\gamma)}_n)(x)$ and $p^{(\gamma)}(x)$
be a sequence of approximate solutions satisfying conditions {\rm (A.1)}--{\rm (A.3)} and {\rm (H)},
and
$$
e_i(\gamma)\rightarrow 0 \qquad \mbox{as $\gamma\rightarrow \infty$}
$$
in the distributional sense for $i=1,2$. Then there exists a
subsequence (still denoted by) $(u^{(\gamma)}, p^{(\gamma)})(x)$ that converges {\it a.e.}
to a weak solution $(\bar{u}, \bar{p})$
of the homogeneous incompressible Euler equations as $\gamma\rightarrow \infty$ :
\begin{equation}\label{HIE}
\begin{cases}
\mbox{\rm div}\,\bar{u}=0,\\
\mbox{\rm div}(\bar{u}\otimes\bar{u})+\nabla\bar{p}=0.
\end{cases}
\end{equation}
\end{proposition}
	
\noindent
\textbf{Proof}. From Theorem $\ref{thm2.1}$,
we know that $(u^{(\gamma)}, p^{(\gamma)})$ converges to $(\bar{u}, \bar{p})$ as $\gamma\rightarrow \infty$.
For the approximate continuity equation, we see that, for any test function $\phi\in C^{\infty}_{c}$,
\begin{eqnarray}
\int e_1(\gamma)\, \phi\, {\rm d}x &=& \int \phi\, \mbox{\rm div}(\rho^{(\gamma)} u^{(\gamma)})\, {\rm d}x\nonumber\\
&=& -\int \nabla\phi\cdot u^{(\gamma)} \rho^{(\gamma)}\, {\rm d}x\nonumber\\
&=& -\int \nabla\phi\cdot u^{(\gamma)} (p^{(\gamma)})^{\frac{1}{\gamma}}\, {\rm d}x.
\end{eqnarray}
Letting $\gamma\rightarrow\infty$, we conclude
\begin{equation}
\int \nabla\phi\cdot \bar{u}\, {\rm d}x=0,
\end{equation}
which implies $(\ref{HIE})_1$ in the distributional sense.
With a similar argument, we can show that $(\ref{HIE})_2$ holds
in the distributional sense.
	
\begin{proposition}[Convergence of approximate solutions for the full Euler flow]\label{thm3.2}
Let $\rho^{(\gamma)}(x)$, $u^{(\gamma)}(x)=(u^{(\gamma)}_1, \cdots, u^{(\gamma)}_n)(x)$, and $p^{(\gamma)}(x)$
be a sequence of approximate solutions satisfying conditions {\rm (A.1)}--{\rm (A.3)} and {\rm (F.1)}--{\rm (F.2)},
and
\begin{eqnarray*}
&&e_i(\gamma)\rightarrow 0 \qquad \mbox{for $i=1,2$},\\
&&(p^{(\gamma)})^{-1}\Big(e_3(\gamma)-u^{(\gamma)}\cdot e_2(\gamma)
+\frac{|u^{(\gamma)}|^2}{2}e_1(\gamma)\Big)\rightarrow 0
\end{eqnarray*}
in the distributional sense as $\gamma\rightarrow \infty$.
Then there exists a subsequence (still denoted by)  $(\rho^{(\gamma)}, u^{(\gamma)}, p^{(\gamma)})(x)$
that converges a.e. to a weak solution $(\bar{\rho}, \bar{u}, \bar{p})$ of
the inhomogeneous incompressible Euler equations as $\gamma\rightarrow \infty$:
\begin{equation}\label{IHICE}
			\begin{cases}
				\mbox{\rm div}\,\bar{u}=0,\\[1mm]
				\mbox{\rm div}(\bar{\rho}\bar{u})=0,\\[1mm]
				\mbox{\rm div}(\bar{\rho}\bar{u}\otimes\bar{u})+\nabla\bar{p}=0.
			\end{cases}
\end{equation}
\end{proposition}
	
\noindent\textbf{Proof}.
From a direct calculation, we have
	\begin{equation}
		\mbox{div}\big((p^{(\gamma)})^{\frac{1}{\gamma}}u^{(\gamma)}\big)=\frac{\gamma-1}{\gamma}(p^{(\gamma)})^{\frac{1}{\gamma}-1}
		\Big(e_3(\gamma)-\sum_{i=1}^n u^{(\gamma)}\cdot e_2(\gamma)+\frac{|u^{(\gamma)}|^2}{2}e_1(\gamma)\Big).
	\end{equation}
	
Then, for any test function $\phi\in C^{\infty}_{c}$, we find
\begin{eqnarray}
&&\int \frac{\gamma-1}{\gamma}(p^{(\gamma)})^{\frac{1}{\gamma}-1}
		\Big(e_3(\gamma)
		-\sum_{i=1}^n u^{(\gamma)}\cdot e_2(\gamma)+\frac{|u^{(\gamma)}|^2}{2}e_1(\gamma)\Big) \phi\, {\rm d}x \nonumber\\
		&=& -\int \nabla\phi\cdot u^{(\gamma)} (p^{(\gamma)})^{\frac{1}{\gamma}}\, {\rm d}x.
\end{eqnarray}
Taking $\gamma\rightarrow\infty$, we have
\begin{equation}
		\int \nabla\phi\cdot \bar{u}\, {\rm d}x=0,
\end{equation}
which implies $(\ref{IHICE})_1$ in the distributional sense.

The fact that $(\ref{IHICE})_2$ and $(\ref{IHICE})_3$ hold in the distributional sense
can be shown similarly from $e_j(\gamma)\rightarrow 0$ as $\gamma\rightarrow \infty$, $j=1,2$, respectively.
	
\smallskip
\begin{remark}\label{rem6}
The main difference between Propositions
$\ref{thm2.2}$ and $\ref{thm3.2}$
is that,
when $\gamma\rightarrow\infty$,
the compressible
homentropic Euler equations converge to the homogeneous incompressible
Euler equations with the unknown variables $(u, p)$,
while
the full Euler equations converge to the inhomogeneous
incompressible Euler equations with the unknown variables $(\rho, u, p)$.
Furthermore, the incompressibility of the limit for the homentropic case
follows directly from the approximate continuity equation $(\ref{3.1})_1$,
while the incompressibility for the full Euler case is from a combination
of all the equations in \eqref{3.2}.
\end{remark}
	
There are various ways to construct approximate solutions by either
numerical methods or analytical methods such as numerical/analytical
vanishing viscosity methods.
As direct applications of the compactness framework,
we now present two examples in \S 3--\S 4 for establishing
the incompressible limit for the multidimensional steady
compressible Euler flows
through infinitely long nozzles.
	
\section{Incompressible Limit for Two-Dimensional Steady Full Euler Flows
in an Infinitely Long Nozzle}
	
In this section, as a direct application of the compactness framework
established in Theorem \ref{thm3.1}, we establish the incompressible
limit of steady subsonic full Euler flows
in a two-dimensional, infinitely long nozzle.
	
The infinitely long nozzle is defined as
	\begin{equation*}
		\Omega=\{(x_1,x_2)\, :\, f_1(x_1)<x_2<f_2(x_1), \,-\infty<x_1<\infty\},
	\end{equation*}
with the nozzle walls $\partial\Omega:=W_1\cup W_2$, where
	\begin{equation*}
		W_i=\{(x_1,x_2)\, :\, x_2=f_i(x_1)\in C^{2,\alpha}, ~-\infty<x_1<\infty\}, \qquad i=1,2.
	\end{equation*}
	
\begin{figure}[htbp]
		\small \centering
		\includegraphics[width=10cm]{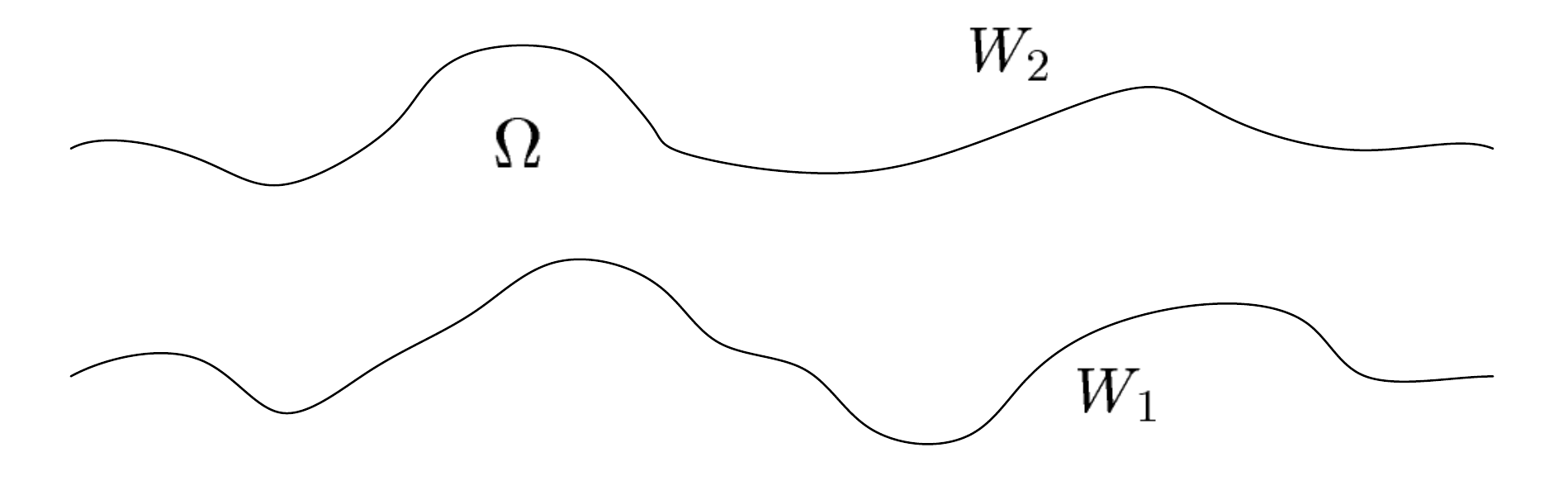}
		\caption{Two-dimensional infinitely long nozzle}
		\label{Fig01}
	\end{figure}

Suppose that $W_1$ and $W_2$ satisfy
\begin{align}\label{cdx-3}
&f_2(x_1)>f_1(x_1) ~ \qquad \mbox{for} ~x_1\in(-\infty, \infty),\nonumber\\[1mm]
&f_1(x_1)\rightarrow 0, \quad f_2(x_1)\rightarrow 1 \qquad \mbox{as} ~x_1\rightarrow -\infty, \nonumber\\[1mm]
&f_1(x_1)\rightarrow a, \quad f_2(x_1)\rightarrow b>a
		\qquad \mbox{as} ~x_1\rightarrow \infty,
\end{align}
and there exists $\alpha>0$ such that
\begin{equation}\label{cdx-4}
		\|f_i\|_{C^{2,\alpha}(\mathbb{R})}\leq C,
		\qquad i=1,2,
\end{equation}
for some positive constant \emph{C}.
It follows that $\Omega$ satisfies the uniform exterior sphere condition
with some uniform radius $r>0$. See Fig \ref{Fig01}.
	
	\bigskip
	
Suppose that the nozzle has impermeable solid walls so that the flow satisfies
the slip boundary condition:
\begin{equation}\label{cdx-5}
u \cdot\nu =0 \qquad \mbox{on} ~\partial \Omega,
\end{equation}
where
$\nu$
is
the unit outward normal to the nozzle wall.
	
It follows from $(\ref{CFE})_1$ and $(\ref{cdx-5})$ that
\begin{equation}\label{cdx-6}
		\int_{\Sigma}\, (\rho u) \cdot l \, {\rm d}s\equiv m
\end{equation}
holds for some constant $m$, which is the mass flux,
where $\Sigma$ is any curve transversal to the $x_1$--direction,
and {$l$} is the normal of $\Sigma$ in the positive $x_1$--axis
direction.
	
We assume that the upstream entropy function is given, {\it i.e.},
\begin{equation}\label{cdx-8}
		\frac{\rho}{p^{1/\gamma}}\longrightarrow S_-(x_2)
		\qquad \mbox{as} ~x_1\rightarrow -\infty,
\end{equation}
and the upstream Bernoulli function is given, {\it i.e.},
\begin{equation}\label{cdx-11}
		\frac{|u|^2}{2}+\frac{\gamma p}{(\gamma-1)\rho}\longrightarrow B_-(x_2)
		\qquad\mbox{as} ~x_1\rightarrow -\infty,
\end{equation}
where $S_-(x_2)$ and $B_-(x_2)$ are the functions defined on $[0,1]$.
	
\medskip
$\mathbf{Problem~1}(m, \gamma)$:
{\it	
Solve the full Euler system $(\ref{CFE})$ with the boundary condition $(\ref{cdx-5})$,
	the mass flux condition $(\ref{cdx-6})$,
	and the asymptotic conditions $(\ref{cdx-8})$--$(\ref{cdx-11})$.
}

\medskip	
Set
\begin{equation*}
		\underline{S}=\inf_{x_2\in[0,1]}S_-(x_2),
		\qquad \underline{B}=\inf_{x_2\in[0,1]}B_-(x_2).
\end{equation*}
	
For this problem, the following theorem has been established in
Chen-Deng-Xiang \cite{ChenDX}.
	
\begin{theorem}\label{thm5.1}
Let the nozzle walls $\partial \Omega$ satisfy
$(\ref{cdx-3})$--$(\ref{cdx-4})$, and let
$\underline{S}>0$ and $\underline{B}>0$.
Then there exists $\delta_0>0$ such that,
if $\|(S_--\underline{S},B_--\underline{B})\|_{C^{1,1}([0,1])}\leq \delta$
for $0<\delta\leq\delta_0$, $(S_-B_-)'(0)\leq 0$, and $(S_-B_-)'(1)\geq 0$,
there exists $\hat{m}\geq 2\delta_0^{\frac{1}{8}}$ such that,
for any $m\in(\delta^{\frac{1}{4}},\hat{m})$,
there is a global solution ({\it i.e.} a full Euler flow)
$(\rho, u, p) \in C^{1,\alpha}(\overline{\Omega})$
of $\mathbf{Problem~1}(m, \gamma)$ such that the following hold:
		
\medskip
{\rm (i)} Subsonic state and horizontal direction of the velocity: The flow is uniformly
subsonic with positive horizontal velocity in the whole nozzle, {\it i.e.},
\begin{equation}\label{cdx3}
\sup_{\bar{\Omega}}(|u|^2-c^2)<0, \quad u_1>0  \qquad \mbox{in} ~\bar{\Omega};
\end{equation}
		
{\rm (ii)} The flow satisfies the following asymptotic behavior in the far field:
As $x_1\rightarrow -\infty$,
\begin{equation}\label{cdx4}
			p\rightarrow p_->0, \qquad
			u_1\rightarrow u_-(x_2)>0, \qquad
			(u_2,\rho)\rightarrow (0,\rho_-(x_2;p_-)),
\end{equation}
\begin{equation}\label{cdx5}
			\nabla p\rightarrow 0, \qquad
			\nabla u_1\rightarrow (0, u_-'(x_2)),
			\qquad \nabla u_2\rightarrow 0,
			\qquad \nabla\rho\rightarrow (0,\rho_-'(x_2;p_-))
\end{equation}
uniformly for $x_2\in K_1\Subset(0,1)$, where
$\rho_-(x_2;p_-)=p_-^{\frac{1}{\gamma}} S_- (x_2)$,
the constant $p_-$ and function $u_-(x_2)$ can be determined by $m$, $S_-(x_2)$,
and $B_-(x_2)$ uniquely.
\end{theorem}
	
Next, we take the incompressible limit of the full Euler flows.
	
\begin{theorem}[Incompressible limit of two-dimensional full Euler flows]\label{thm5.2}
Let $(\rho^{(\gamma)}, u^{(\gamma)}, p^{(\gamma)})(x)$ be the corresponding
sequence of solutions to $\mathbf{Problem~1} (m^{(\gamma)}, \gamma)$.
Then, as $\gamma\rightarrow \infty$, the solution
sequence possesses a subsequence (still denoted by)
$(\rho^{(\gamma)}, u^{(\gamma)}, p^{(\gamma)})(x)$
that converges strongly {\it a.e.} in $\Omega$ to
a vector function $(\bar{\rho}, \bar{u}, \bar{p})(x)$ which is a weak solution
of $(\ref{ICIHE})$.
Furthermore, the limit solution $(\bar{\rho}, \bar{u}, \bar{p})(x)$ also satisfies
the boundary condition $(\ref{cdx-5})$ as the normal trace of the divergence-measure
field $u$ on the boundary in the sense of Chen-Frid {\rm \cite{Chen7}}.
\end{theorem}
	
\noindent{\bf Proof}.  We divide the proof into four steps.
	
\medskip
1. From $(\ref{CFE})$, we can obtain the following linear transport parts:
\begin{eqnarray}\label{3.2.1}
\begin{cases}
			\partial_{x_1}\big((p^{(\gamma)})^{\frac{1}{\gamma}} u_1^{(\gamma)}\big)
			+\partial_{x_2}\big((p^{(\gamma)})^{\frac{1}{\gamma}} u_2^{(\gamma)}\big)=0,\\[2mm]
			\partial_{x_1}\big((p^{(\gamma)})^{\frac{1}{\gamma}}  B^{(\gamma)}u_1^{(\gamma)} \big)
			+\partial_{x_2}\big((p^{(\gamma)})^{\frac{1}{\gamma}}  B^{(\gamma)}u_2^{(\gamma)}\big)=0,\\[2mm]
			\partial_{x_1}\big((p^{(\gamma)})^{\frac{1}{\gamma}} S^{(\gamma)} u_1^{(\gamma)}\big)
			+\partial_{x_2}\big((p^{(\gamma)})^{\frac{1}{\gamma}} S^{(\gamma)} u_2^{(\gamma)}\big)=0.
\end{cases}
\end{eqnarray}

From $(\ref{3.2.1})_1$, we can introduce the potential function $\psi^{(\gamma)}$:
	\begin{eqnarray}
		\begin{cases}
			\partial_{x_1}\psi^{(\gamma)}=-(p^{(\gamma)})^{\frac{1}{\gamma}} u_2^{(\gamma)},\\[2mm]
			\partial_{x_2}\psi^{(\gamma)}=(p^{(\gamma)})^{\frac{1}{\gamma}} u_1^{(\gamma)}.
		\end{cases}
	\end{eqnarray}
	
From the far-field behavior of the Euler flows,
we can define
$$
\psi^{(\gamma)}_-(x_2):=\lim\limits_{x_1\rightarrow-\infty}\psi^{(\gamma)}(x_1, x_2).
$$
Since both the upstream Bernoulli and entropy functions are given,
$B^{(\gamma)}$ and $S^{(\gamma)}$ have the following expressions:
$$
	B^{(\gamma)}(x)=B_{-}((\psi^{(\gamma)}_{-})^{-1}(\psi^{(\gamma)}(x))),
	\qquad S^{(\gamma)}(x)=S_{-}((\psi^{(\gamma)}_{-})^{-1}(\psi^{(\gamma)}(x))),
$$
where $(\psi^{(\gamma)}_{-})^{-1}\psi^{(\gamma)}(x)$ is a function
from $\Omega$ to $[0,1]$, and
\begin{eqnarray}
		\begin{cases}
			B_-=\frac{u_-^2}{2}+\frac{\gamma}{\gamma-1}\frac{p_-}{\rho_-},\nonumber\\[2mm]
			S_-=\rho_-p_-^{-\frac{1}{\gamma}},
		\end{cases}
\end{eqnarray}
with uniformly upper and lower bounds with respect to $\gamma$.
	
Since the flow is subsonic so that the Mach number $M^{(\gamma)}\le 1$, then we have
\begin{equation}
|u^{(\gamma)}|<\sqrt{\frac{2(\gamma-1)}{\gamma+1}\max B_-},
\end{equation}
and
\begin{equation}
		\left(\frac{2(\gamma-1)}{\gamma(\gamma+1)}\min (B_- S_-)\right)^{\frac{\gamma}{\gamma-1}}
		< p^{(\gamma)}
		\le\left(\frac{\gamma-1}{\gamma}\max (B_-S_-)\right)^{\frac{\gamma}{\gamma-1}}.
\end{equation}
Since $|u^{(\gamma)}|^2$ and $p^{(\gamma)}$ are uniformly bounded,
we conclude that $|u^{(\gamma)}|^2$ and $p^{(\gamma)}$ are uniformly bounded in $L^1_{loc}(\Omega)$. 	
Thus, conditions {\rm (A.1)}--{\rm (A.2)} are satisfied.
It is observed that, even though the lower bound of pressure $p^{(\gamma)}$ may tend to zero as $\gamma\to \infty$ {with polynomial rate}, so that (F.1) holds for any bounded domain.
	
\medskip
2. For fixed $x_1$, $(\psi^\gamma_{-})^{-1}(\psi^\gamma(\cdot))$
can be regarded as a backward characteristic map with
$$
\frac{\partial ((\psi^{(\gamma)}_{-})^{-1}(\psi^{(\gamma)}))}{\partial x_2}
=\frac{(p^{(\gamma)})^{\frac{1}{\gamma}} u^{(\gamma)}_1}{p_-^{\frac{1}{\gamma}}u_{-}}>0.
$$
The uniform boundedness and positivity of $p_-^{\frac{1}{\gamma}}u_{-}$
and $(p^{(\gamma)})^{\frac{1}{\gamma}} u^{(\gamma)}_1$
implies that
the map is not degenerate.
Then we have
\begin{eqnarray}
\begin{cases}
			\partial_{x_1}S^{(\gamma)}(x)
			=-S_{-}'((\psi^{(\gamma)}_{-})^{-1}(\psi^{(\gamma)}(x)))
			\frac{(p^{(\gamma)})^{\frac{1}{\gamma}} u^{(\gamma)}_2}{p_-^{\frac{1}{\gamma}}u_{-}},\\[2mm]
			\partial_{x_2}S^{(\gamma)}(x)
			=S_{-}'((\psi^{(\gamma)}_{-})^{-1}(\psi^{(\gamma)}(x)))
			\frac{(p^{(\gamma)})^{\frac{1}{\gamma}} u^{(\gamma)}_1}{p_-^{\frac{1}{\gamma}}u_{-}}.
\end{cases}
\end{eqnarray}
Thus, $S^{(\gamma)}$ is uniformly bounded in $BV$, which implies its strong convergence.
Then condition {\rm (F.2)} follows.
	
\medskip
3. Similar to \cite{ChenHuangWang}, the vorticity sequence $\omega^{(\gamma)}:=\partial_{x_1}u_2^{(\gamma)}-\partial_{x_2}u_1^{(\gamma)}$ can be written as
\begin{eqnarray}
\begin{cases}
\partial_{x_1} B^{(\gamma)}=  u_2^{(\gamma)} \omega^{(\gamma)}-\frac{\gamma}{\gamma-1}(\rho^{(\gamma)})^{-2}(p^{(\gamma)})^{\frac{\gamma+1}{\gamma}}\partial_{x_1} S^{(\gamma)},\\[2mm]
\partial_{x_2} B^{(\gamma)}=  - u_1^{(\gamma)} \omega^{(\gamma)}-\frac{\gamma}{\gamma-1}(\rho^{(\gamma)})^{-2}(p^{(\gamma)})^{\frac{\gamma+1}{\gamma}}\partial_{x_2} S^{(\gamma)}.
\end{cases}
\end{eqnarray}
By direct calculation, we have
\begin{eqnarray}
\omega^{(\gamma)}&=&\frac{1}{|u^{(\gamma)}|^2}
\Big(u_2^{(\gamma)}\big(\partial_{x_1} B^{(\gamma)}+\frac{\gamma}{\gamma-1}(\rho^{(\gamma)})^{-2}(p^{(\gamma)})^{\frac{\gamma+1}{\gamma}} \partial_{x_1} S^{(\gamma)}\big)\nonumber\\
&&\qquad\qquad -u_1^{(\gamma)}\big(\partial_{x_2} B^{(\gamma)}+\frac{\gamma}{\gamma-1}(\rho^{(\gamma)})^{-2}(p^{(\gamma)})^{\frac{\gamma+1}{\gamma}}\partial_{x_2} S^{(\gamma)}\big)\Big)\nonumber\\
&=&-\frac{1}{p_-^{\frac{1}{\gamma}} u_-}
\Big((p^{(\gamma)})^{\frac{1}{\gamma}} B_-' +\frac{\gamma}{\gamma-1}(p^{(\gamma)})^{\frac{\gamma+2}{\gamma}} (\rho^{(\gamma)})^{-2} S_-'\Big),
\end{eqnarray}
which implies that $\omega^\varepsilon$ as a measure sequence is uniformly bounded
so that it is compact in  $H^{-1}_{loc}$.
Therefore, the flows satisfy condition {\rm(A.3)}.
	
\smallskip
Then {Proposition $\ref{thm3.2}$} immediately implies that
the solution sequence has a subsequence (still denoted by)
$(\rho^{(\gamma)}, u^{(\gamma)}, p^{(\gamma)})(x)$
that converges {\it a.e.} in $\Omega$ to
a vector function $(\bar{\rho}, \bar{u}, \bar{p})(x)$
as $\gamma\to \infty$.
	
\smallskip
4. Since $\bar{u}$ is uniformly bounded,
the normal trace $\bar{u}\cdot \nu$ on $\partial \Omega$ exists and is in $L^\infty(\partial \Omega)$
in the sense of Chen-Frid \cite{Chen7}.
On the other hand, for any $\phi\in C^{\infty}(\mathbb{R}^2)$, we have
\begin{equation}
\langle (\bar{u}\cdot \nu)|_{\partial \Omega}, \phi\rangle= \int_{\Omega} \bar{u}(x) \cdot \nabla \phi(x)\, {\rm d}x
+ \int_{\Omega} \phi\mbox{div}\,\bar{u} \, {\rm d}x.
\end{equation}
Since
$\int_{\Omega} \phi\,\mbox{div}\,\bar{u} \, {\rm d}x=0$,
and
\begin{equation}
\int_{\Omega}\bar{u}(x) \cdot \nabla \phi(x)\, {\rm d}x
		=\lim\limits_{\gamma \rightarrow \infty}\int_{\Omega}((p^{(\gamma)})^{\frac{1}{\gamma}} u^{(\gamma)})(x) \cdot \nabla \phi(x)\, {\rm d}x=0,
\end{equation}
then we have
\begin{equation}
\langle (\bar{u}\cdot \nu)|_{\partial \Omega}, \phi\rangle=0,
\end{equation}
for any $\phi\in C^{\infty}(\mathbb{R}^2)$.
By approximation, we conclude that
the normal trace $(\bar{u}\cdot \nu)|_{\partial \Omega} = 0$ in $L^\infty(\partial\Omega)$.
This completes the proof.
	
\smallskip
\begin{remark}
In the two-dimensional homentropic case, the subsonic results in \cite{CX,Xin4} can also be
extended to the incompressible limit by using {Proposition $\ref{thm2.2}$}.
\end{remark}
	
\smallskip
\section{Incompressible Limit for the Three-Dimensional Homentropic Euler Flows in an Infinitely Long Axisymmetric Nozzle}
	
\begin{figure}[htbp]
		\small \centering
		\includegraphics[width=6cm]{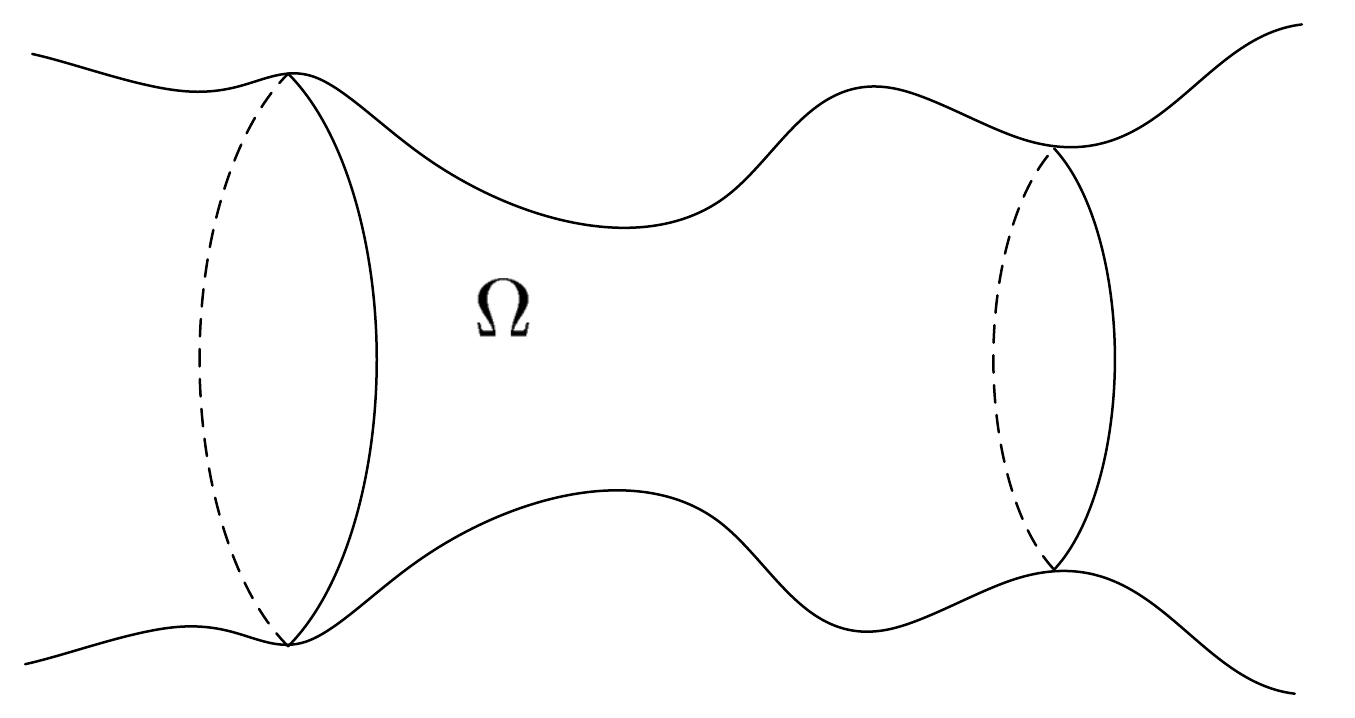}
		\caption{Infinitely long axisymmetric nozzle}
		\label{Fig3}
\end{figure}

\smallskip
We consider Euler flows through an infinitely long axisymmetric nozzle
in $\R^3$ given by
\begin{equation*}
\Omega=\{(x_1,x_2,x_3)\in\mathbb{R}^3\, : \, 0\leq \sqrt{x_2^2+x_3^2}<f(x_1),~
		-\infty<x_1< \infty\},
\end{equation*}
	where $f(x_1)$ satisfies
	\begin{align}\label{dd-2}
		&f(x_1)\rightarrow 1 \qquad \mbox{as} ~x_1\rightarrow-\infty,\nonumber\\[1mm]
		&f(x_1)\rightarrow r_0  \qquad \mbox{as} ~x_1\rightarrow \infty,\nonumber\\[1mm]
		&\|f\|_{C^{2,\alpha}(\mathbb{R})}\leq C \qquad \mbox{for}
		~\mbox{some} ~\alpha>0 ~ \mbox{and} ~C>0, \\[1mm]
		& \inf_{x_1\in\mathbb{R}}f(x_1)=b>0.
	\end{align}
	See  Fig. \ref{Fig3}.
	
\bigskip

\medskip
The boundary condition is set as follows:
Since the nozzle wall is solid,
the flow satisfies the slip boundary condition:
\begin{equation}\label{dL-3}
		u\cdot\nu=0 \qquad \mbox{on} ~\partial \Omega,
\end{equation}
where
$\nu$ is the unit outward
normal to the nozzle wall.
The continuity equation in $(\ref{CFE})_1$ and the
boundary condition $(\ref{dL-3})$ imply that
the mass flux
\begin{equation}\label{dd-4}
\int_\Sigma \, (\rho u)\cdot l\, {\rm d}s\equiv m_0
	\end{equation}
remains for some positive constant $m_0$,
where $\Sigma$ is any surface transversal to the $x_1$--axis direction,
and
$l$ is the normal of $\Sigma$ in the positive $x_1$--axis direction.
	
\smallskip
In Du-Duan \cite{DD}, axisymmetric flows without swirl are considered for
the fluid density $\rho=\rho(x_1,r)$ and
the velocity
$$
u=(u_1, u_2, u_3)=(U(x_1,r), V(x_1, r)\frac{x_2}{r}, V(x_1, r)\frac{x_3}{r})
$$
in the cylindrical coordinates,
where $u_1$, $u_2$, and $u_3$ are the axial velocity, radial velocity,
and swirl velocity, respectively,
and $r=\sqrt{x_2^2+x_3^2}$.
Then, instead of $(\ref{CISE})$, we have
\begin{equation}\label{dd-6}\left\{ \begin{split}
			&\partial_{x_1}(r\rho U)+\partial_r(r\rho V)=0,\\[1mm]
			&\partial_{x_1}(r\rho U^2)+\partial_r(r\rho UV)_r+r\partial_{x_1}p=0,\\[1mm]
			&\partial_{x_1}(r\rho UV)+\partial_r(r\rho V^2)_r+r\partial_rp=0.
		\end{split} \right.
\end{equation}
Rewrite the axisymmetric nozzle as
\begin{equation*}
		\Omega=\{(x_1,r)\, :\, 0\leq r<f(x_1), ~-\infty<x_1<\infty\}
\end{equation*}
with the boundary of the nozzle:
\begin{align*}
\partial \Omega=\{(x_1,r)\, :\, r=f(x), ~-\infty<x_1<\infty\}.
\end{align*}
	
\smallskip
The boundary condition $(\ref{dL-3})$ becomes
\begin{equation}\label{dd-7}
		(U,V,0)\cdot \tilde{\nu}=0 \qquad \mbox{on} ~\partial \Omega,
	\end{equation}
where
$\tilde{\nu}$
is the unit outer normal of the nozzle wall
in the cylindrical coordinates.
The mass flux condition $(\ref{dd-4})$ can
be rewritten in the cylindrical coordinates as
\begin{equation}\label{dd-9}
		\int_\Sigma (r\rho U, r \rho V,0)\cdot \tilde{l}\, {\rm d}s\equiv m:=\frac{m_0}{2\pi},
	\end{equation}
where $\Sigma$ is any curve transversal to the $x_1$-axis direction,
and
$\tilde{l}$
is the unit normal of $\Sigma$.
	
\smallskip
Notice that the quantity
$$
B=\frac{\gamma}{\gamma-1}\rho^{\gamma-1}+\frac{U^2+V^2}{2}
$$
is constant along each streamline.
For the homentropic Euler flows in the axisymmetric nozzle,
we assume that the upstream Bernoulli is given, that is,
\begin{equation}\label{dd-13}
		\frac{\gamma}{\gamma-1}\rho^{\gamma-1}+\frac{U^2+V^2}{2}\longrightarrow B_-(r) \qquad \mbox{as} ~x_1\rightarrow -\infty,
\end{equation}
where $B_-(r)$ is a function defined on $[0, 1]$.
	
\smallskip
Set
\begin{equation}
\underline{B}=\inf\limits_{r\in[0,1]} B_-(r),
\quad \sigma=||B'_-||_{C^{0,1}([0,1])},
\end{equation}
We denote the above problem as $\mathbf{Problem~2} (m, \gamma)$. It is shown in \cite{DD} that
	
\begin{theorem}\label{thm4.5}
Suppose that the nozzle satisfies $(\ref{dd-2})$.
Let the upstream Bernoulli function $B(r)$
satisfy $\underline{B}>0$, $B'(r)\in C^{1,1}([0,1])$, $B'(0)=0$,
		and $B'(r)\geq 0$ on $r\in [0,1]$.
		Then we have
\begin{enumerate}
\item[(i)]
There exists $\delta_0>0$ such that, if $\sigma\leq\delta_0$,
then there is $\hat{m}\leq 2\delta_0^{\frac{1}{8}}$.
For any $m\in (\delta_0^{\frac{1}{4}},\hat{m})$,
there exists a global $C^1$--solution ({\it i.e.} a homentropic Euler flow)
$(\rho, U, V)\in C^1(\overline{\Omega})$ through the nozzle with mass flux
condition $(\ref{dd-9})$ and the upstream asymptotic condition $(\ref{dd-13})$.
Moreover, the flow is uniformly subsonic, and the axial velocity is always
positive, {\it i.e.},
\begin{equation}\label{dd2}
\sup_{\bar{\Omega}}(U^2+V^2-c^2)<0 \quad \mbox{and} \quad U>0
\qquad \mbox{in} ~\overline{\Omega}.
\end{equation}
			
\item[(ii)] The subsonic flow satisfies the following properties: As $x_1\rightarrow-\infty$,
\begin{equation}\label{dd4}
\rho\rightarrow\rho_->0, \quad \nabla\rho\rightarrow 0, \quad p\rightarrow \rho_-^\gamma,\quad (U, V)\rightarrow (U_-(r), 0), \quad \nabla U\rightarrow(0,U_-'(r)),
\end{equation}
uniformly for $r\in K_1\Subset(0,1)$,
where $\rho_-$ is a positive constant, and $\rho_-$ and $U_-(r)$ can
be determined by $m$ and $B(r)$ uniquely.
\end{enumerate}
\end{theorem}
	
As above, we have the following incompressible limit theorem for this case.
	
\begin{theorem}[Incompressible limit of three-dimensional Euler flows through an axisymmetric nozzle]\label{thm4.3}
Let	$u^{(\gamma)}=(u^{(\gamma)}_1, u^{(\gamma)}_2, u^{(\gamma)}_3)$,
and $p^{(\gamma)}=(\rho^{(\gamma)})^{\gamma}$ be the corresponding solutions
to $\mathbf{Problem~2}~(m^{(\gamma)}, \gamma)$.
Then, as $\gamma \rightarrow \infty$, the solution
sequence  possesses a subsequence (still denoted by)
$(u^{(\gamma)}, p^{(\gamma)})$
that converges strongly {\it a.e.} in $\Omega$ to
a vector function
$(\bar{u}, \bar{p})$ with $\bar{u}=(\bar{u}_1, \bar{u}_2, \bar{u}_3)$
which is a weak solution of $(\ref{ICHE})$.
Furthermore, the limit solution $(\bar{u}, \bar{p})$ also satisfies
the boundary conditions $(\ref{dL-3})$
as the normal trace of the divergence-measure field
$(\bar{u}_1, \bar{u}_2, \bar{u}_3)$
on the boundary in the sense of Chen-Frid {\rm \cite{Chen7}}.
\end{theorem}
	
\noindent{\bf Proof}.
For the approximate solutions, $B^{(\gamma)}$
satisfy
\begin{eqnarray}
\partial_{x_1}(r U^{(\gamma)} (p^{(\gamma)})^{\frac{1}{\gamma}}  B^{(\gamma)})
+\partial_r (r V^{(\gamma)} (p^{(\gamma)})^{\frac{1}{\gamma}}  B^{(\gamma)}) = 0.
	\end{eqnarray}
Based on the equation:
$$
\partial_{x_1}(r U^{(\gamma)} (p^{(\gamma)})^{\frac{1}{\gamma}})+\partial_r (r V^{(\gamma)} (p^{(\gamma)})^{\frac{1}{\gamma}}) = 0,
$$
we introduce $\psi^{(\gamma)} $ as
\begin{eqnarray}
\begin{cases}
\partial_{x_1}\psi^{(\gamma)} =-r V^{(\gamma)} (p^{(\gamma)})^{\frac{1}{\gamma}},\\[2mm]
\partial_{r}\psi^{(\gamma)} =r U^{(\gamma)} (p^{(\gamma)})^{\frac{1}{\gamma}}.
\end{cases}
\end{eqnarray}
From the far-field behavior of the Euler flows,
we define
$$
\psi^{(\gamma)}_-(r):=\lim\limits_{x_1\rightarrow-\infty}\psi^{(\gamma)}(x_1, r).
$$

Similar to the argument in Theorem \ref{thm5.2}, $(\psi^{(\gamma)}_{-})^{-1}(\psi^{(\gamma)})$ are
nondegenerate maps.
A direct calculation yields
\begin{eqnarray*}
B^{(\gamma)}(x_1, x_2, x_3)=B_{-}((\psi^{(\gamma)}_{-})^{-1}(\psi^{(\gamma)}(x_1, \sqrt{x_2^2+x_3^2}))),
\end{eqnarray*}
with
\begin{equation}
B^{(\gamma)}_-=\frac{U_-^2}{2}+\frac{\gamma}{\gamma-1}\rho_-^{\gamma-1},
\end{equation}
	
Similar to the previous case, the flow is subsonic so that the Mach number $M^{(\gamma)}\le 1$,
\begin{equation}
|(U^{(\gamma)},V^{(\gamma)})|<\sqrt{2\max B_-},
\end{equation}
and
\begin{equation}
\left(\frac{2(\gamma-1)}{\gamma(\gamma+1)}\min B_- \right)^{\frac{\gamma}{\gamma-1}}
< p^{(\gamma)}
\le\left(\frac{\gamma-1}{\gamma}\max B_-\right)^{\frac{\gamma}{\gamma-1}}.
\end{equation}
From \eqref{total-energy}, we have
\begin{equation}
\frac{2}{\gamma(\gamma+1)}\min B_-
< E^{(\gamma)}
\le\frac{(\gamma-1)\gamma+2}{2\gamma}\max B_-.
\end{equation}
Therefore, conditions {\rm (A.1)}--{\rm (A.2)} and (H)
are satisfied for any bounded domain.

\smallskip
On the other hand, the vorticity $\omega^{(\gamma)}$ has the following expressions:
\begin{eqnarray}
\begin{cases}
\omega_{1,2}^{(\gamma)}=\partial_{x_1}u_2^{(\gamma)}-\partial_{x_2}u_1^{(\gamma)}
=\frac{x_2}{r}(\partial_{x_1}V^{(\gamma)}-\partial_{r}U^{(\gamma)}),\\[2mm]
\omega_{2,3}^{(\gamma)}=\partial_{x_2}u_3^{(\gamma)}-\partial_{x_3}u_2^{(\gamma)}=0,\\[2mm]
\omega_{3,1}^{(\gamma)}=\partial_{x_3}u_1^{(\gamma)}-\partial_{x_1}u_3^{(\gamma)}
=-\frac{x_3}{r}(\partial_{x_1}V^{(\gamma)}-\partial_{r}U^{(\gamma)}).
\end{cases}
\end{eqnarray}
A direct calculation yields
\begin{equation}
\partial_{x_1}V^{(\gamma)}-\partial_{r}U^{(\gamma)}
=-\frac{r}{p_-^{\frac{1}{\gamma}} U_-}(p^{(\gamma)})^{\frac{1}{\gamma}} B_-',
\end{equation}	
which implies that $\omega^{(\gamma)}$ is uniformly bounded in the bounded measure space and (A.3) is satisfied.
	
Then the sequence $(u^{(\gamma)}, p^{(\gamma)})(x)$ satisfies
conditions  {\rm (A.1)}--{\rm (A.3)} and {\rm (H)}.
Moreover, $(\ref{CISE})$ holds for 	$(u^{(\gamma)}, p^{(\gamma)})(x)$.
	
Similar to Theorem \ref{thm5.2}, we conclude that there exists
a subsequence (still denoted by) $(u^{(\gamma)}, p^{(\gamma)})$ that
converges to a vector function $(\bar{u},\bar{p})$ {\it a.e.} in $\Omega$
satisfying \eqref{ICHE}
in the distributional sense.
	
\smallskip
Since $\bar{u}$ is uniformly bounded,
the normal trace $\bar{u}\cdot \nu$ on $\partial \Omega$ exists and is in $L^\infty(\partial \Omega)$
in the sense of Chen-Frid \cite{Chen7}.
On the other hand, for any $\phi\in C^{\infty}(\mathbb{R}^2)$, we have
\begin{equation}
\langle (\bar{u}\cdot \nu)|_{\partial \Omega}, \phi\rangle= \int_{\Omega} \bar{u}(x) \cdot \nabla \phi(x)\, {\rm d}x
+ \int_{\Omega} \phi\,\mbox{div}\,\bar{u} \, {\rm d}x.
\end{equation}
Since
$\int_{\Omega} \phi\,\mbox{div}\,\bar{u} \, {\rm d}x=0$,
and
\begin{equation}
\int_{\Omega}\bar{u}(x) \cdot \nabla \phi(x)\, {\rm d}x=0,
\end{equation}
then we have
\begin{equation}
\langle (\bar{u}\cdot \nu)|_{\partial \Omega}, \phi\rangle=0,
\end{equation}
for any $\phi\in C^{\infty}(\mathbb{R}^2)$.
By approximation, we conclude that
the normal trace $(\bar{u}\cdot \nu)|_{\partial \Omega} = 0$ in $L^\infty(\partial\Omega)$.
This completes the proof.

\begin{remark}
For the full Euler flow case, the subsonic results of \cite{Duan-Luo} can be also extended to
the incompressible limit by Proposition $\ref{thm3.2}$.
\end{remark}

\medskip
\noindent {\bf Acknowledgments:}
The research of Gui-Qiang G. Chen was supported in part by
the UK EPSRC Science and Innovation
Award to the Oxford Centre for Nonlinear PDE (EP/E035027/1),
the UK EPSRC Award to the EPSRC Centre for Doctoral Training
in PDEs (EP/L015811/1), and
the Royal Society--Wolfson Research Merit Award (UK).
The research of Feimin Huang was supported in part by
NSFC Grant No. 10825102 for distinguished youth scholars,
and the National Basic Research Program of China (973 Program)
under Grant No. 2011CB808002.
The research of Tianyi Wang was supported in part
by the China Scholarship Council  No. 201204910256
as an exchange graduate student at the University of Oxford,
the UK EPSRC Science and Innovation Award to the Oxford Centre
for Nonlinear PDE (EP/E035027/1),
and the NSFC Grant No. 11371064;
He would like to thank Professor Zhouping Xin for the helpful discussions.
Wei Xiang was supported in part by the UK EPSRC Science and Innovation Award
to the Oxford Centre for Nonlinear PDE (EP/E035027/1),
the CityU Start-Up Grant for New Faculty 7200429(MA),
and the General Research Fund of Hong Kong
under GRF/ECS Grant 9048045 (CityU 21305215).

\end{document}